
%
%

\magnification 1200
\font\abs=cmr9
\font\ccc=cmcsc10

\def\tens{\otimes}
\def\fraz#1#2{{\strut\displaystyle #1\over\displaystyle #2}}

\def\ii#1{\item{$\phantom{1}#1\,.~$}}
\def\jj#1{\item{$#1\,.~$}}
\def\c{{\bf C}}
\def\R{{\bf R}}
\def\Z{{\bf Z}}
\def\N{{\bf N}}
\def\d{\Delta}

\def\pro #1#2 {{\ccc (#1.#2) Proposition.}\phantom{X}}
\def\dfn #1#2 {{\ccc (#1.#2) Definition.}\phantom{X}}
\def\cor #1#2 {{\ccc (#1.#2) Corollary.}\phantom{X}}
\def\lem #1#2 {{\ccc (#1.#2) Lemma.}\phantom{X}}
\def\rem #1#2 {{\ccc (#1.#2) Remark.}\phantom{X}}
\def\thm #1#2 {{\ccc (#1.#2) Theorem.}\phantom{X}}
\def\v#1{| #1 \rangle}
\def\dim {{\sl Proof.}\phantom{X}}
\def\fidi{\hskip5pt \vrule height4pt width4pt depth0pt}
\def \n {\bar n}
\def \v {\bar v}
\def \z {\bar z}
\def \vn { \{ v,n\} }
\def \vnn { \{ v,\n \} }
\def \nn { \{ n,\n \} }
\def \dv {\partial _v }
\def \dn {\partial _n }
\def \dnn {\partial _{\n} }

\def \o {\omega}
\def \ed {{\cal E}(2)}
\def \eo {{\cal E}_{\o}(2)}
\def \co {{\cal C}_{\o}}
\def \I {{\cal I}}
\def \A {{\cal A}}

\hsize= 15 truecm
\vsize= 22 truecm
\hoffset= 0.2 truecm
\voffset= 0.0 truecm
\baselineskip= 14 pt
\footline={\hss\tenrm\folio\hss} \pageno=1

\vglue 3truecm
\baselineskip= 18 pt
\centerline {\bf QUANTUM PLANES AND QUANTUM CYLINDERS}
\centerline{\bf FROM POISSON HOMOGENEOUS SPACES.}
\baselineskip= 14 pt
\bigskip
\bigskip
\centerline{{\it N. Ciccoli}}
\medskip
\centerline{\abs {Faculteit der Wiskunde en Informatica,\footnote
\dag {on leave from Dipartimento di Matematica, Universit\'a di Bologna.
Italy.} }}
\centerline {\abs {Universiteit van Amsterdam.}}
\centerline {\abs {The Netherlands.}}
\bigskip
\bigskip
\baselineskip= 12 pt
{\bf Abstract.}
{\abs Quantum planes and a new quantum cylinder are obtained
as quantization of Poisson homogeneous spaces of two different Poisson
structures on classical euclidean group E(2).
}
\baselineskip= 14 pt
\bigskip
\bigskip
\noindent {\bf 1. Introduction.}
\bigskip
\smallskip
The concept of homogeneous space of a group is maybe one of the most
widespread mathematical concept, lying, for example, at the very
foundation of harmonic analysis and of symmetries of physical systems.
{}From the very beginning of the theory of quantum groups, clarifying
the concept of a quantum homogeneous space for a quantum group has been
considered of the uttermost importance (see [11] for example), although up
to now a complete theory is still lacking.
The purpose of this work is to study homogeneous quantum spaces of
euclidean quantum groups
through the analysis of the "semiclassical" limit
Poisson structure on the classical group.
Our aim is to verify to what extent the
results in [12] and [13] relating covariant Poisson structures on the sphere
and the one parameter family of quantum spheres [11] are
still valid in this case (although we will not deal with the analytical
aspects as it is extensively done in those references).
There are two different versions of the Euclidean quantum group, introduced
in [15] and [3] respectively; in this second reference it is also  shown how
they can be  both obtained
through contraction procedure from $SU_h (2)$. The way in which they are
related is explained in [1]. Quite a lot of work has been
done on these groups, whose interest lies both in being an easy example of the
nonsemisimple and noncompact case and in physical applications. For
example in [16]
the standard euclidean quantum group is treated at an analytical level, in
[15] its relations with $q$-special functions are first studied and later
in [2] and [9] (last of a series of paper and useful for further references)
$q$-harmonic analysis is more thouroghly investigated.
The roots of unity theory has been
dealt with in [5] for the standard case, and has been claimed and proved
to be trivial in the non standard case respectively in [1] and [4].
The quantum homogeneous spaces have been studied
in the standard case in [2], through duality arguments involving the
quantized universal enveloping algebra.
The semiclassical limit of these results is here shown to fit in the
framework of the classification of covariant Poisson structures. This
is the content of section 2. Furthermore
in the nonstandard case the semiclassical limit suggests quite
naturally that homogeneous spaces should be of cylinder-type. The covariant
Poisson structures on this cylinder are classified in section 3 and in
section 4 an explicit quantum cylinder is given through generators and
relations, verifying conditions explained in [6] for quantum homogeneous
spaces.

\bigskip
\noindent {\bf 2. The standard Poisson algebraic $\ed$ and quantum planes.}
\bigskip
The classical 2-dimensional Euclidean group is usually written in
complex notation as the space of matrices of the form:

$$ \left( \matrix {v&n\cr 0&1\cr }\right) $$
\bigskip
Its function (polynomial) algebra can thus be seen as the commutative
algebra on generators $v$, $\v$, $n$, $\n$, with the additional
relation $v \v =1$. The matrix form gives us immediatly also the
Hopf algebra structure that we explicitly write down as follows:
$$ \d v=v\tens v \qquad \qquad \d \v =\v \tens \v$$
$$ \d n =\v \tens n+n\tens 1 \qquad
\d \n =v\tens \n +\n \tens 1$$
$$ S(v)=\v \qquad S(n)=-vn\qquad S(\n )=-\v \n$$
$$ \varepsilon (v)=1 \qquad \varepsilon (n)=\varepsilon (\n ) =0$$
and the usual $*$-structure is given by:
$$ v^* =\v \qquad n^* =\n \qquad *^2=Id$$.

\smallskip
Let us define on this function algebra the following
quadratic Poisson bracket (standard bracket in what follows):
$$\{ v,n\} =vn \qquad \{ v,\n \} = v\n$$
$$\{ n,\n \} =n\n \qquad \{ v,\v \} =0$$
\rem 21 These formulas define a Poisson bracket on the algebra of
polynomial functions on $\ed$.  That's why we talk about an {\it algebraic}
Poisson structure and not of a Lie-Poisson structure that should be given
on the whole algebra of smooth functions. In the following we will some times
drop the adjective algebraic.

\medskip
\pro 22 {\it The formulas above define a Poisson algebraic structure on $\ed$.}
\smallskip
\dim
Let us write the Poisson bivector  corresponding to the bracket as:
$$ w(v,n,\n )=vn\dv \wedge \dn +v\n \dv \wedge \dnn
+n\n \dn \wedge \dnn $$
By direct computation this Poisson bivector verifies
the  multiplicativity property:
$$w(gg')=((L_g )'_{g'}\tens (L_g)'_{g'})w(g')+((R_{g'})'_g \tens
(R_{g'})'_g)(w(g))$$
where $L_g$ and $R_g$ stands respectively for left and right
translations in $\ed$.
\fidi
\medskip
Differentiating the Poisson bivector at the origin gives a coalgebra structure
on the Lie algebra $e(2)$. Choosing $J=\dv$ , $X=\dn$ , $Y=\dnn$ as generators
of the Lie algebra (obviously considering them around the identity) we have the
coalgebra structure:
$$ \delta (J )=0 \qquad \delta (X )=J \wedge X$$
$$ \delta (Y )=J \wedge Y$$
It is easy to show that such structure is non-coboundary.

The Poisson bivector above mentioned defines the associated
homomorphism $B_w$ from the cotangent
to the tangent bundle of $\ed$ that can be easily calculated to be:
$$dv \mapsto -v(n\dn +\n \dnn )= X_v$$
$$ dn \mapsto n(-v\dv +\n \dnn )=X_n$$
$$ d\n \mapsto \n (n\dn +v\dv )=X_{\n}$$
Thus we have the differentiable distribution of tangent subspaces:
$$(v,n,\n )\rightarrow (X_v ,X_n ,X_{\n} )$$
whose integral manifolds are the symplectic leaves of the Poisson structure.
\smallskip
\rem 23 {\it Geometric description of symplectic leaves.}
Let us observe that for a generic point of the space $(v, n,\n )$ we
have :
$$ vn\n X_v +\n X_n +nX_{\n } =0$$
and thus the tangent space is generically 2-dimensional.

All the points of the form $(v,0,0)$ are $0-dim$ symplectic
leaves of this Poisson distribution. The set of all this point is
nothing but ${\bf S}^1$ viewed as a subgroup of $\ed$, and thus,
being an union of leaves, it is a Poisson algebraic subgroup. It is
not difficult
to see that this subgroup, together with all its finite subgroup
are the only Poisson algebraic subgroup of $\ed$.
\medskip
 At this point it could be interesting to calculate the primitive
spectrum of standard $Fun_q (\ed)$ to see if it is in a good relation
with the symplectic foliation. The primitive spectrum of this
algebra is not very difficult to explicit after we've noticed
that it can be given the form of
an iterated skew (twisted) polynomial ring (see [7], [8]).
\medskip
\cor 24 {\it The $\ed$-homogeneous space ${\R }^2$ has a natural
Poisson structure induced by the algebraic Poisson structure on $\ed$.}
\medskip
The projection map $\pi :\ed \rightarrow \R ^2$ is nothing but the map
that sends $(v,n,\n )\rightarrow (\fraz 1 2 (n+\n) ,-i\fraz 1 2 (n-\n) )$.
Calling $z$ and $\z$  the coordinates on the plane it is then immediate
that the Poisson structure on the plane is:
$$ \{ z,\z \} =z\z $$
By definition such a structure is $(\ed ,w)$-covariant. We recall that
by covariance we mean that the action map $\phi :\ed \times \R ^2\rightarrow
\R ^2$ is a Poisson map when the first space is given the product
Poisson structure.
\smallskip
\rem 25 The Poisson structure just defined on the plane has a symplectic
foliation consisting in two families of $0-dim$ leaves parametrized by
points of two orthogonal lines and four $2-dim$ leaves
separated by those points.
Observing that the usual quantum plane structure $F_q$, the algebra on
two $q$-commuting generators, is a quantization of this Poisson bracket,
it's no surprise that the primitive spectrum of this algebra is in
bijective correspondence with the foliation (see [14] for the explicit
expression of the primitive spectrum).
\smallskip
The next problem is that of classifying all $(\ed ,w)$-covariant Poisson
structures on the plane. We will follow [17] where the covariance condition
is rewritten at the infinitesimal level. We will denote with:
$$ \phi :\ed \times \R ^2 \rightarrow \R ^2 $$
the action map, with $\phi _x =\phi (\cdot ,x)\quad \forall x\in R^2 $ and
$\phi _g = \phi(g,\cdot )\quad \forall g\in \ed$. With these notations
condition 2.2
of [17] states that covariant structures on the plane are in 1-1 correspondence
with elements $\rho \, \in \bigwedge ^2 T_0 (\R ^2)$ such that

$$(\phi _0)_{*,e} \delta (X) +X\cdot \rho =0
\qquad \forall X \in {\cal G}_{0} $$

where ${\cal G}_0$ is the tangent algebra of the rotation subgroup -
the stabilizer of $0$ in $\R ^2$. In fact such a $\rho$ can be
extended to a Poisson bivector field on $\R ^2$ simply as:

$$\rho (x)=(\phi _0) _{*,g}(w(g)) + (\phi _g )_{*,0}(\rho )$$

where $g\in \ed$ is such that $g\cdot 0=x$
. Now, as ${\cal G}_0$ is nothing but the algebra generated by
$J$ and $\delta (J)=0$, the above condition rewrites as $J\cdot \rho =0$
(simply an invariance condition) which is always verified. We can
then take $\rho = k\partial _1 \wedge \partial _2$, where the derivatives
are calculated in $0$. The corresponding bivector field is:
$$  (a,b)\mapsto k\partial _1 \wedge \partial _2$$
\medskip
and the Poisson bracket on the plane is:
$$ \{ z,\z \} =z\z +k. $$
Note that when $\rho =0$ we obtain exactly the Poisson structure on the plane
given by Corollary 2.4.
We have then just proved the following
\smallskip
\pro 26 {\it There is a one-parameter family of covariant Poisson
structures on the plane, with respect to the standard Poisson structure
on $\ed$, given by:
$$ \{ z,\z \} =z\z +k.$$
}
\medskip
\rem 27 In the case in which $k\ne 0$ the symplectic foliation changes
drastically. The 0-dimensional symplectic leaves are the points of
the hyperbola $z\z =-k$, and they divide the plane in three $2$-dimensional
symplectic leaf. Following this geometric picture we will refer to this
case as the case of hyperbolic covariant structures.

\bigskip
\bigskip
\noindent {\bf 3. The nonstandard Poisson structure and the cylinder.}

\bigskip
We still deal with the same function algebra but we want to consider another
family of Poisson structure (to which we'll refer as non standard Poisson
structures) given in [10] and whose quantization
gives the so-called non standard
euclidean quantum group.
$$\vn =\o (1-v)$$
$$\vnn =-\o (v^2 -v)$$
$$\nn =\o (n-\bar n )$$
where $\o$ is a nonzero complex number.

That this bracket gives a Poisson algebraic structure on $\ed $ (or better a
family of isomorphic structures) is proven in [10].

The infinitesimal counterpart of the bracket is the coproduct on the Lie
algebra $e(2)$ given by:

$$ \delta (P_1 )=0 \qquad \delta (P_2 )=\o P_2 \wedge P_1 \qquad
\delta (J)=\o J\wedge P_2 .$$

This is a coboundary coproduct by taking:

$$ \delta (X)=ad_X r, \qquad r=\o J\wedge P_2 .$$

We will call $w'$ the Poisson bivector and
$B_{w'}$ the associated homomorphism
between tangent and cotangent spaces is:

$$dv\mapsto \o (v-1)\dn +\o (v^2 -v)\dnn =\o (v-1)(\dn +v\dnn )=X_v$$
$$dn\mapsto \o (1-v)\dv +\o (\n -n)\dnn =X_n$$
$$d\n \mapsto \o (v-v^2 )\dv +\o (n-\n )\dn =X_{\n}$$

from which we have the distribution of tangent subspaces that integrates
to simplectic leaves.
\medskip

\medskip
\rem 31 For every $ (v, n, \n )$ the relation:
$$ (v-v^2 )X_n +(v-1)X_{\n } +(\n -n)X_v =0$$
holds, showing that the distribution is at most two-dimensional at every
point. If we restrict to point of ${\bf S}^1$ given by $n=\n =0$ as
in the previous case, we see that the distribution is
exactly two-dimensional in these points
and thus every point of the form $(v,0,0)$ is contained in a two-dimensional
symplectic leaf. Thus ${\bf S}^1$ is not a Poisson subgroup of non standard
Poisson $\ed$.
If we restrict to points $v=1$ and $n=\n$ we see that the distribution
vanish and thus all these points are $0-dim$ symplectic leaves. Again
we're dealing with a subgroup (that in the matrix form corresponds to
upper triangular unimodular $2\times 2$ real matrices) and it is isomorphic
to $\R$. Thus $\R$ is a Poisson subgroup for  the non standard
Poisson structure on $\ed$. This implies that the corresponding
homogeneous space, a cylinder, has a canonical covariant Poisson structure.
One further remark could be that in this case also the discrete infinite
groups are Poisson subgroup. However we have to notice that, as the
Poisson structure is defined only at the algebraic level, such subgroups
are not closed, thus can not be given by annhilating an element of the function
algebra.

Again it may be interesting to confront the whole symplectic foliation
with the primitive ideal structure of its quantization, to compute
which one could use the fact that $\ed$ is an iterated differential
polynomial ring (as noted in [4]).
\medskip
\pro 32 {\it The cylinder $\cal C$, viewed as a homogeneous space of the
non standard Poisson Euclidean group, can be parametrized
as having the function algebra generated
by $v$, $\v$ and $\v \n -vn  = m$ with the covariant Poisson structure given
by:
$$ \{ v,m \} =-\o (v^2 -1)$$ . }
\smallskip
\dim
It is enough to show that the two functions $v$ and $m$ parametrize
the cylinder, the rest will follow immediatly. That the function $v$
is invariant with respect to the action of $ \cal A$ is obvious, and an
easy calculation shows also the invariance of $m$.
\fidi
\medskip
\rem 33 From the explicit Poisson bracket of $3.3$ one can obtain the
corresponding symplectic foliation on the cylinder. The distribution of
tangent subspaces is:
$$ dv\mapsto \o (v^2 -1) \partial _m =X_v$$
$$ dm \mapsto -\o (v^2 -1) \dv =X_m$$
and thus all the points on the line $v=1$ are $0-dim$ symplectic leaves
and all the other points belong to a unique $2-dim$ leaf.
\medskip
\rem 34 Now we want to give a classification of covariant structures
on the cylinder  as in $2.6$. Again multiplicative $(\ed ,w')$-Poisson
structures are characterized by their value in one point, say $v=1, m=0$
. In this point a bivector has the form $\varrho = k\dv \wedge \partial _m$
and, with notations as in $2.5$,  has to fulfill the invariance condition:
$$ (\phi _{(1,0)})_* \delta (P_1 )+P_1 \cdot \varrho = 0$$
that reduces to the condition $P_1 \cdot \varrho = 0$  due to triviality
of $\delta (P_1 )$ . This condition is
always verified. This implies that the space of covariant
Poisson structures has the form of an affine 1-dimensional
space given adjoining to the structure of Proposition 3.2 any invariant
bivector field extending $\varrho $. Explicitely the possible Poisson
structures are listed in the following:
\medskip
\pro 35 {\it There is a one parameter family of covariant Poisson
structures on the cylinder, with respect to the non standard Poisson
$\ed$, given by:
$$ \{ v,m\} = -\o (v^2 -1) +k.$$}
\medskip
\rem 36 As for the plane the case $k\ne 0$ shows a completely
different symplectic foliation. Just remark that the bracket is trivial
when:
$$ \o v^2 +\o +k=0$$
from which we obtain the solutions:
$$ v=\pm \sqrt (1+\o ^{-1} k)$$
The condition that $v$ is on the unit circle can be expressed as $v\v =1$.
Except from the singular
situation in which $\o ^{-1
} k \in i\R$,  there is
always a special value of $k\ne 0$ for which $1+\o ^{-1}k$ (and thus its
square roots) belong to the unit circle. Thus for every complex $\o$
which is not purely imaginary there are always two values of $v$ for which
the Poisson bracket is degenerate. The corresponding symplectic foliation
is then given by two lines of points and two $2-dim$ symplectic leaves
between them. For all the other values of $k$ the distribution is two
dimensional in every point and so the Poisson bracket induces a
symplectic form on the whole cylinder.

\bigskip

\bigskip
\noindent {\bf 4. Quantum homogeneous cylinder}
\bigskip
We want now to pass from the semiclassical situation of Poisson homogeneous
spaces to quantum homogeneous spaces, as defined for example in [6]. The
standard $\eo$ has been treated in [2], where connections with the
theory of spherical functions are also explained, so that in what follows
we will restrict to
the nonstandard case where the quantization of the given Poisson structure
can be obtained, as in [10], simply substituting the Poisson bracket with
the commutator, leaving the coalgebra structure unchanged.
Let us just remark the connections between [2] and Proposition $2.6$.
The natural covariant Poisson structure is exactly the semiclassical
limit of the usual quantum plane, and is, in fact the only quantum
homogeneous spaces obtained by "quotient" with respect to a proper
quantum subgroup, incidentally the quantization of the unique Poisson
subgroup of $\ed$. The quantum hyperboloid of [2] has as semiclassical
limit the covariant Poisson structure of $2.6$ with $k=-2$ (and the
analysis of symplectic leaves give a clear meaning to the appearing
of the world "hyperbolic" in its name). The fact that covariant Poisson
structures constitutes a complete one-parameter family either suggests
the possibility of other quantum planes of which the two given in [13]
should be, in a sense, paradigmatic, or asks for an explanation of the
failure of quantizing Poisson homogeneous spaces for values of
$k$ different from $0$ and $-2$.

Let's move to the nonstandard case. Let us recall ([1]) that the nonstandard
quantum group is the Hopf algebra with the same coproduct, counit and
antipode as the standard quantum group and commutation relations obtained
substituting the Poisson bracket with the Lie bracket in formulas at the
beginning of chapter 3.

\medskip
The classical projection $\pi :G\rightarrow M$ from one Poisson-Lie group
to its homogeneous space correspond to an injective map on the function
algebra level ${\hat \pi}:{\cal F}(M)\rightarrow {\cal F}(G)$. The remark
just stated above says that we can see this map as a map between the
quantizations and requiring it to be an algebra map we find the quantized
commutation relations for the Poisson homogeneous space. It is obvious how
it works in the quantum plane case. This implies that we can give the
following:
\medskip
\dfn 41 {\it The quantum cylinder $\co$ is the algebra generated
by the elements $v$, $\v$ and $m$ with the relations:
$$ v\v =\v v=1 \qquad vm=mv-\o (v^2 -1) \qquad \v m=m\v +\o (\v -\v^2 )$$
and with the $*$-structure:
$$ v^* =\v \qquad m^* =-m$$.}
\medskip
\pro 42 {\it The quantum cylinder is a quantum homogeneous space of
 nonstandard $\eo$ in the sense of [6] i.e. it is a  $*$-invariant right
coideal
. Furthermore the elements $\{ v^r m^s : r\in \Z , s\in \N \}$ provide
a basis as a vector space for this algebra.}
\smallskip
\dim
Obviously the algebra $\co$ is $*$-invariant.
To prove that it is a
right coideal in $Fun_q (\ed )$ is enough to perform straightforward
calculations on the generators.
We have:
$$ \d v=v\tens v  \in  \eo \tens \co $$
$$ \d m =1\tens m+\v \n \tens \v -vn \tens v  \in \eo \tens
\co .$$
To prove the linear independance of elements $v^r m^s$ let us first observe
that $\co$ is a skew polynomial ring ([7]), and precisely $\co = R[m;\delta ]$
where $R=\c [v,v^{-1}]$ and $\delta (v)=(v^2 -1)$. Thus it is possible to
define the degree in $m$, that we will denote $deg_m$, for every element
in $\co$ and prove that if $n$ and $u$ are in $\co$ then:
$$ deg_m (nu)=deg_m (n)+deg _m (u).$$
It is then clear that there cannot be linear relations involving
monomials with $deg _m \ne 0$. On the other side linear relations involving
only elements of $m$-degree $0$ are linear relations in $\c [v,v^{-1}]$ and
thus they should be trivial.

Let us observe that $\co$ being a skew polynomial ring is easily shown
to verify good algebraic properties through very general arguments.
For example it is a Noetherian integral domain.
\fidi
\medskip
In [6] it is shown how, in a fixed Hopf algebra,  it is possible
to construct a galoisian reciprocity
between $*$-invariant subalgebras and right coideals on one side
and Hopf $*$-ideals on the other side through the assignments:
$$ \Sigma :B\mapsto \A _B =\langle (S^n -\varepsilon {\bf 1})
   (b), b\in B, n\in \Z \rangle$$
$$ \Pi :\A \mapsto B_{\A}=\{ b: (\pi \tens id)\circ \d (b)=1\tens b\}$$
where $\pi$ is the projection from the whole Hopf algebra onto the
quotient with respect to $\A$.
It is then natural then to give the following:
\medskip
\dfn 43 {\it The closure of a quantum homogeneous space $\cal C$
in the Hopf algebra ${\cal H}_q$ is the homogeneous space
${\cal B}=\Pi \circ \Sigma ({\cal C})$.}

\medskip

\pro 44 {\it The closure of the quantum cylinder $\co$ is the quantum
homogeneous
space corresponding to the quantum subgroup of non standard $\eo$ given
by the Hopf-$*$-ideal
$\, \I = \langle v-1, n-\n \rangle $.}
\smallskip
\dim
The notation of last line of the proposition means that $\I$ is the ideal
generated by elements enclosed in brackets.
We want to verify that $\I$ is really a Hopf-$*$-ideal.
First we have to verify that this ideal is a bilateral coideal i.e.:
$$ \d \I = \I \tens \eo +\eo \tens \I $$
$$ \varepsilon (\I )=0$$
It is enough to prove it for the generators. The second condition is
trivial. The first follows from:
$$ \d (v-1)=v\tens v-1\tens 1 =(v-1)\tens 1+v\tens (v-1)$$
$$ \d (n-\n )=\v \tens n+n\tens 1-v\tens \n -\n \tens 1
=\v \tens (n-\n )+(n-\n )\tens 1+(\v -v)\tens \n$$
We next have to prove $S$-invariance. Again on generators we have:
$$ S(v-1)=(\v -1)=-\v (1-v); \qquad S(n-\n )=m=\v (\n -n)-(v-1)(\v +1)n$$
so that $S(\I )\subset \I$.
As for the $*$-structure we have:
$$(v-1)^* =\v -1=-\v (v-1)  \in \I$$
$$ (n-\n )^* =n^* -\n ^* =-(n-\n )  \in \I.$$

Next we have to verify that the closure of $\co$, as defined in $(4.3)$,
coincides
with the set:
$$ B_{\I} =\{ b\in \eo : (\pi \tens id)\circ \d (b)=1\tens b \}$$
where $\pi :\eo \rightarrow \eo
/ \I$ is the
natural projection.
Direct calculations show that indeed $v$ and $m$ belongs to this set and thus
that $\co$ is contained in $B_{\I}$. For example:
$$ (\pi \tens id) (\d m)=1\tens m+\pi (\v (\n -n))\tens \v
   +\pi ((\v -v)n)\tens v.$$
The stabilizer ideal of $\co$ is:
$$ \A _{\co} = \Sigma (\co) =\langle S^n (c)-\varepsilon (c)1, b\in \co,
\, n\in \Z \rangle $$
Let us note that
$(S-\varepsilon {\bf 1})(m)=n-\n $ and $(S-\varepsilon {\bf 1})(\v)=v-1$
so that the generators
of $\I$ belongs to $\A _{\co}$ proving $\I \subset \A_{\co}$. On the other
hand from general arguments $\co \subset B_{\I}$ implies the opposite
inclusion so that:
$$ \A_{\co }=\I .$$
This proves that the closure of $\co$
is exactly $B_{\I}$.
\fidi

\bigskip
\bigskip
\noindent {\bf 5. Conclusions}
\bigskip
Analizing the Poisson structures and their Poisson homogeneous spaces
for the euclidean group we have recovered the semiclassical limit of
results in [2] about quantum homogeneous planes, showing that,
in this limit there is a one-parameter family of planes, as it is
the case for spheres ([12], [13]). Furthermore we have determined a
quantization of a Poisson structure on the cylinder compatible
with a quantum homogeneous structure with respect to nonstandard
quantum euclidean group, showing that it can be recovered
from an explicit quantum subgroup. Again we have showed a one parameter family
of covariant Poisson structures on this space.
The method of looking at the semiclassical limit situation to have hints
towards the general quantum homogeneous spaces could perhaps give
indication in all those situations in which the function algebra is
explicitely given through generators and relations.

\bigskip
\noindent {\bf Acknowledgements}
\smallskip
I am indebted with Dr. B. Drabant and Prof. R. Giachetti for useful
discussions on the subject and to Prof. T.H. Koornwinder for careful
reading of an early version of the manuscript.

\bigskip
\bigskip

\centerline{\bf References.}
\bigskip

\ii 1 Ballesteros A., Celeghini E., Giachetti R., Sorace E. and
    Tarlini M., {\it An R-matrix approach to the quantization of the
    Euclidean group E(2).}  J. Phys. A: Math. Gen., {\bf 26}, 7495 (1993).
\ii 2 Bonechi F., Ciccoli N., Giachetti R., Sorace E., Tarlini M.,
    {\it Free $q$-Schr\"odinger equation from homogeneous spaces of the
    $2-dim$ Euclidean quantum group.}  Comm. Math. Phys. {\it in press.}
\ii 3 Celeghini E., Giachetti R., Sorace E. and Tarlini M., {\it $3-dim$
    quantum groups from contractions of $SU_q (2)$.} J. Math. Phys.
    {\bf 31} 2548 (1990) and  {\it Contractions of Quantum Groups.}
    Lecture notes in Mathematics $1510$, {\bf 221} (Berlin, Springer, 1992).
\ii 4 Ciccoli N., {\it Integer forms of nonsemisimple quantum groups},
    preprint, (1995).
\ii 5 Ciccoli N., Giachetti R., {\it The euclidean quantum algebra at
      roots of unity.} Lett. Math. Phys., {\it in press.}
\ii 6 Dijkhuizen M.S. and Koornwinder T.H., {\it Quantum homogeneous spaces,
    duality and quantum 2-spheres.} Geom. Ded. {\bf 52}
    291 (1994).
\ii 7 Goodearl K. R., Letzter E. S., {\it Mem. Amer. Math. Soc.},
    {\bf 521}, (1994).
\ii 8 Jordan D. A., {\it Iterated skew polynomial rings and quantum groups.}
      Journ. of Alg., {\bf 156}, 194 (1993).
\ii 9 Koelink H. T., {\it The quantum group of plane motion and the Hahn-Exton
      $q$-Bessel function.} Duke Math. Journ., {\bf 76}, 483 (1994).
\jj {10} Maslanka P., {\it The $E_q (2)$ group via direct quantization of the
    Lie-Poisson structure and its Lie algebra.} Journ. Math. Phys.
    {\bf 35}, 1976 (1994).
\jj {11} Podles P., {\it Quantum spheres.} Lett. Math. Phys.
    {\bf 14}, 193 (1987).
\jj {12} Sheu A. J. L., {\it Quantization of the Poisson $SU(2)$ and its
     Poisson homogeneous space, the $2$-sphere.} Comm. Math. Phys.
    {\bf 135}, 217 (1991).
\jj {13} Sheu A. J. L., {\it Quantum Poisson $SU(2)$ and quantum Poisson
    spheres.} Contemp. Math. {\bf 134}, 247 (1992).
\jj {14} Smith. A. P., {\it Quantum groups: an introduction and survey for
    ring theorists.} in ``Noncommutative rings'', {\it MSRI Publications},
    (Montgomery S. and Small L. W. eds.) {\bf 24}, 131 (Berlin, Springer,
1992).
\jj {15} Vaksman L.L. and Korogodski L.I., {\it An algebra of bounded functions
     on the quantum group of the motions of the plane and $q$-analogues of
     Bessel functions.}  Sov. Math. Dokl. {\bf 39},
    173 (1989).
\jj {16} Woronowicz S., L., {\it Unbounded elements affiliated with $C^*$-
    algebras and non compact quantum groups.} Comm. Math. Phys., {\bf 136},
    399 (1991) and {\it Quantum $SU(2)$ and $E(2)$ groups: contraction
procedure    .} Comm. Math. Phys., {\bf 149}, 637 (1992).
\jj {17} Zakrzewski S., {\it Poisson homogeneous spaces. } Proceedings of
    the Karpacz winter school on
    theoretical physics 1994, {\it to appear}. Also, {\it hep-th 9412101}.

\bye